\documentclass[12pt]{amsart}

\usepackage{amssymb,latexsym,amssymb,amsthm}

\newcommand{\R}{{\mathcal R}}
\newtheorem{theorem1}{Theorem}
\newtheorem{theorem2}[theorem1]{Theorem}
\newtheorem*{remark}{Remark}
\newtheorem{lemma}{Lemma}
\newtheorem{lemma2}[lemma]{Lemma}
\newtheorem{lemma3}[lemma]{Lemma}
\title{The argument of the Riemann $\xi$ function off the critical line}
\author{Xiannan Li}
\email{xli@math.stanford.edu}
\address{Department of Mathematics, Stanford University, Stanford, CA 94305}

\begin{document}
\maketitle
\begin{abstract}
We examine the behaviour of the zeros of the real and imaginary parts of $\xi(s)$ on the vertical line $\Re s = 1/2+\lambda$, for $\lambda \neq 0$.  This can be rephrased in terms of studying the zeros of families of entire functions $A(s) = \frac{1}{2} (\xi(s+\lambda) + \xi(s - \lambda))$ and $B(s) = \frac{1}{2i} (\xi(s+\lambda) - \xi(s - \lambda))$.  We will prove some unconditional analogues of results appearing in \cite{Lag}, specifically that the normalized spacings of the zeros of these functions converges to a limiting distribution consisting of equal spacings of length $1$, in contrast to the expected GUE distribution for the same zeros at $\lambda = 0$.  We will also show that, outside of a small exceptional set, the zeros of $\Re \xi(s)$ and $\Im \xi(s)$ interlace on $\Re s = 1/2+\lambda$.  These results will depend on showing that away from the critical line, $\arg \xi(s)$ is well behaved.
\end{abstract}

\section{Introduction}
Let $F(s) = \pi^{-s/2}\Gamma(s/2)\zeta(s)$, where as usual, $s = \sigma + it$ is a complex variable, and $\zeta(s)$ is the Riemann zeta function.  Let $\xi(s) = \frac{(s)(s-1)}{2}F(s)$.  Then $F(s)$ and $\xi(s)$ are real on the $\Re s = 1/2$ line.  We are interested in studying the zeros of the real or imaginary parts of $F(s)$ or $\xi(s)$ off the $1/2$ line.  As will become clear in the proofs, there is essentially no difference in the behaviour of $F(s)$ or $\xi(s)$ in this situation.

Let $n_\lambda(T)$ denote the number of zeros of $\Re F(1/2+\lambda + it)$ with $0<t\leq T$.  For convenience, we can also assume that $\lambda > 0$ since $F(s) = F(1-s)$.  Let $N_\lambda(T)$ be the number of zeros of the zeta function in the region $\Re(s) \geq 1/2 + \lambda$ and $0<\Im(s) \leq T$.  Ki \cite{Ki} improved a result of Levinson \cite{Lev} and showed that 
$$n_\lambda(T) = \frac{T}{2\pi}\log \frac{T}{2\pi} - \frac{T}{2\pi} + O((N_\lambda(T)+1) \log T).
$$By an improvement of a result of Selberg by Jutila (see for instance 9.29 of \cite{Titch}), we have that $N_\lambda(T) \ll T^{1- a \lambda}\log T$ where $a$ is any constant less than 1, so that the formula above is indeed an asymptotic formula.  In the course of the proof, Ki also noted that the number of zeros of $\Re F(1/2+\lambda + it)$ in the interval $[T, T+1]$ is $\ll \log T$.  Similar results hold for the zeros of $\Im F(1/2+\lambda + it)$.

The similarities between the results above and basic results on zeros of the zeta function are striking.  However, in \cite{Lag}, Lagarias proves assuming RH for $0<\lambda < 1/2$ and unconditionally for $\lambda \geq 1/2$ that the behaviour of these zeros is quite different from that of the zeros of the zeta function.  Actually, Lagarias proved his results for slightly different functions, in essence for $R(t) = \Re \xi(1/2+\lambda + it)$ and $I(t) = \Im \xi(1/2+\lambda + it)$.  Lagarias proved, unconditionally for $\lambda > 1/2$ and conditionally assuming R.H. for $\lambda < 1/2$, that the distributions of the normalized spaces of zeros of $R(t)$ or $I(t)$ up to height $T$ converge to a limiting distribution which consists of equal spaces of size $1$.  This is in sharp contrast to the conjectured GUE distribution of zeros of $F(s)$ at $\Re s = 1/2$.  Lagarias also showed that the zeros of $I(t)$ and $R(t)$ interlace and are simple.  That is, for any two consecutive zeros of $R(t)$, $\gamma_1$ and $\gamma_2$, there exists exactly one zero of $I(t)$ in $(\gamma_1, \gamma_2)$ and vice versa.

Lagarius's perspective was somewhat different from that of Levinson in that instead of dealing with the real part or imaginery part of certain functions, he instead viewed his result in terms of differences of L-functions.  Specifically he examined the entire functions $A(s) = \frac{1}{2} (\xi(s+\lambda) + \xi(s - \lambda))$ and $B(s) = \frac{1}{2i} (\xi(s+\lambda) - \xi(s - \lambda))$, which can be considered as a family of entire functions depending on $\lambda$.  By the functional equation and the fact that $\xi(\bar{s}) = \bar{\xi(s)}$,
\begin{eqnarray*}
A(1/2+it) &=& \frac{1}{2} (\xi(1/2+it+\lambda) + \xi(1/2+it - \lambda))\\
&=& \frac{1}{2} (\xi(1/2+it+\lambda) + \xi(1/2-it + \lambda))\\
&=& R(t).
\end{eqnarray*}
Similarly, $B(1/2+it) = I(t)$.  Thus the correspondence between these two viewpoints is that the functions agree on the critical line.  

Our results are unconditionally analogues of the results in \cite{Lag}, which depended on RH in certain ranges.  We will prove our results for $A(s)$ and $B(s)$, noting that the proof applies with trivial changes to $\Re F(1/2+\lambda+it)$ and $\Im F(1/2+\lambda+it)$ as well.  It will also provide alternative routes to the results in \cite{Lag}, which we indicate in a later remark.  Our first result is

\begin{theorem1}
Let $A(s) = \frac{1}{2} (\xi(s+\lambda) + \xi(s - \lambda))$ and $B(s) = \frac{1}{2i} (\xi(s+\lambda) - \xi(s - \lambda))$.  Then for $\lambda > 0$, the following holds unconditionally

1.  $A(s)$ and $B(s)$ have almost all of their zeros on the line $\Re s= 1/2$, with an exceptional set of zeros up to height $T$ not on the line of at most $\ll T^{1-a\lambda}\log ^2T$ for any $a<1$.

2.  Consider consecutive imaginary parts of the zeros of $A(s)$ or $B(s)$ to be $0 \leq\gamma_1 \leq \gamma_2 \leq...\leq \gamma_N$ and then let $\tilde{\gamma_i} = \frac{2\pi \gamma_i}{\log (\gamma_i+1)}$ be the normalized zeros.  We may similarly define normalized zeros for those zeros with negative imaginary part.  The distribution of the spacings of normalized zeros of $A(s)$ and $B(s)$ on the critical line, up to height $T$, converges to a limiting distribution consisting of equal spacings of length $1$.  
\end{theorem1}

The first part of the theorem is essentially a result of Ki in \cite{Ki}, and we will sketch a proof of it in the last section.  The second part of the theorem follows easily from the next theorem, which is our main focus.  

\begin{theorem2}
Let $R(t) = \Re \xi(1/2+\lambda + it)$ and $I(t) = \Im \xi(1/2+\lambda + it)$.  For any $B>0$, there exists $E \subset [T, 2T]$ with measure $\ll\frac{T}{\log^BT} $ and containing $\ll \frac{T}{\log^BT}$ zeros of $I(t)$ and $R(t)$ such that outside of this set, the zeros of $I(t)$ and $R(t)$ are simple, interlace and are regularly spaced.  
\end{theorem2}

The result above stated for $t\in [T, 2T]$ implies the result for $t\in [0, T]$, by diadic summation.  This result will be proven in the next section.  Throughout this paper, we shall assume a result that was essentially proven by Ki in \cite{Ki}, namely that the number of zeros in $[T, T+1]$ of $R(t)$ or $I(t)$ is $\ll \log T$.  Ki proved this for the imaginery and real part of $F(s)$, and his proof, which is a standard application of the argument principle, applies without change to our case.  

The result above is not in the strongest possible form.  The exceptional set as stated is of the size $T$ over any large power of $\log T$, but the result can be proven with the power of $\log T$ replaced by a small power of $T$, depending on $\lambda$.  Note that we expect it to be very difficult to prove the above result with no exceptional set at all, since that implies RH.  Our results extend in general to Dirichlet L-functions.
$\\$
\paragraph{Acknowledgements:} I am very grateful to Professor Soundararajan for his guidance throughout the making of this paper.  I would like to thank Bob Hough, Vorrapan Chandee and Leo Goldmakher for editorial comments.  I am also indebted to the referee for helpful suggestions.

\section{Proof of Theorem 2}
For the rest of this note, let $\sigma = 1/2 + \lambda$ for $\lambda > 0$, and $s = \sigma + it$, for $t\in [T, 2T]$.  We start by noting that $R(t)$ has a zero precisely when $\arg \xi (s) \equiv \pi/2$ mod $\pi$.  Similarly $I(t)$ has a zero when $\arg \xi (s) \equiv 0$ mod $\pi$.  Now by Stirling's formula,
\begin{eqnarray*}
&&\log \left( \frac{s(s-1)}{2}\pi^{-\frac{s}{2}}\Gamma\left(\frac{s}{2}\right)\right) \\
&=& -\frac{s}{2}\log \pi + \left(\frac{s}{2} - \frac{1}{2}\right)\log \frac{s}{2} - \frac{s}{2} + \frac{\log 2\pi}{2}  + \log\frac{s(s-1)}{2} + O(1/s),
\end{eqnarray*}
from which we get
$$\arg \left(\frac{s(s-1)}{2}\pi^{-s/2}\Gamma(s/2)\right) = \frac{t}{2}\log \frac{t}{2\pi} - \frac{t}{2} + O(1).
$$
To understand $\arg \xi(s)$, we argue that the above dominates $\arg \zeta(s)$ outside of some exceptional set, so we need to examine $\arg \zeta(s)$.  Say that $t\in [T, 2T]$ for some large $T$.  From Selberg, we have the standard identity (see e.g. 14.20 of \cite{Titch})
\begin{eqnarray*}
\frac{\zeta '}{\zeta}(s) &=& -\sum_{n<x^2}\frac{\Lambda_x(n)}{n^s} + \frac{x^{2(1-s)} - x^{1-s}}{\log x (1-s)^2} + \frac{1}{\log x} \sum_{q\geq 1} \frac{x^{-2q-s} - x^{-2(2q+s)}}{(2q+s)^2} \\
&+& \frac{1}{\log x} \sum_\rho \frac{x^{\rho - s} - x^{2(\rho-s)}}{(s-\rho)^2}\\
&=& -\sum_{n<x^2}\frac{\Lambda_x(n)}{n^s} + \frac{1}{\log x} \sum_\rho \frac{x^{\rho - s} - x^{2(\rho-s)}}{(s-\rho)^2} + O\left(\frac{x}{T^2}  \right),
\end{eqnarray*}
where
$$\Lambda_x(n) = \begin{cases}
\Lambda(n) & \textup{ for }1\leq n\leq x \\
\Lambda(n) \frac{\log\frac{x^2}{n}}{\log x} &\textup{ for } x\leq n< x^2.
\end{cases}
$$
Here, the latter sum is a sum over nontrivial zeros of $\zeta$.  This sum is only bounded nicely if we stay away from zeros of $\zeta$.  Specifically, for $A$ a parameter to be chosen later, define $E_1 = E_1(A)$ to be the union of all rectangles with vertices $\sigma_0 + i\gamma - i\log ^A T$, $\sigma_0 + i\gamma + i\log ^A T$, $1 + i\gamma + i\log ^A T$, and $1+ i\gamma - i\log ^A T$, over all $\gamma \in [T, 2T]$ such that $\zeta(\beta + i\gamma) = 0$ for some $\beta \geq \sigma_0$.  Then $\mu(E_1) \ll N_{\lambda/2}(T) \log^AT \ll \frac{T}{\log^B T}$, and $E_1$ contains $\ll \frac{T}{\log^B T}$ zeros in total, for any $B>0$.  Indeed, $E_1$ consists of $\ll N_{\lambda/2}(T)$ intervals each of which contains $\ll\log^{A+1}T$ zeros of $R(t)$ and $I(t)$, so $E_1$ contains $\ll \frac{T}{\log^B T}$ zeros in total, for any $B>0$.  

Let $\delta$ be a constant satisfying $0 < \delta<1$; we will discuss choices for $\delta$ later.   Set $x = \log^{\frac{2(1+\delta)}{\lambda}} T$.  Then we have the following Lemma.  

\begin{lemma}
There is a choice of $A$ depending only on $\delta$ and $\lambda$ such that outside of $E_1=E_1(A)$,  \begin{equation}
\frac{\zeta '}{\zeta}(s) = -\sum_{n<x^2}\frac{\Lambda_x(n)}{n^s} + O\left(\frac{1}{\log^\delta T}\right).
\end{equation}

\end{lemma}
\begin{proof}
Let $\sigma_0 = 1/2 + \lambda/2$.  Outside of $E_1$, $ \sum_\rho \frac{x^{\rho - s} - x^{2(\rho-s)}}{(s-\rho)^2}$ is small.  As usual, write $\rho = \beta + i\gamma$.  We split the sum over $\beta \geq \sigma_0$ and $\beta < \sigma_0$.  The sum over $\beta \geq \sigma_0$ is 
$$\ll \frac{x^{2(1-\sigma)}}{\log ^{A-1}T}.
$$Here we bounded the numerator by $x^{2(1-\sigma)}$, noted that $|s-\rho| \geq \log^AT$ for $\beta \geq \sigma_0$ and $s \notin E_1$, and used the bound of $\log T$ on the number of zeros $\rho$ of $\zeta$ satisfying $\Im\rho \in [T, T+1]$.  

Similarly, the sum over $\beta < \sigma_0$ is
$$\ll x^{-\lambda/2}\log T \ll \frac{1}{\log^\delta T},
$$by the definition of $x$.  Since $x$ is a power of $\log T$, we may pick $A$ sufficiently large such that the first bound is $\ll \frac{1}{\log^\delta T}$.  Specifically, we set $A = (1-2\lambda)\frac{2(1+\delta)}{\lambda} + 1 + \delta$.  The result then follows.
\end{proof}

To grasp the size of the sum over primes appearing in Equation 1, we will refer to Corollary 9.5 in \cite{IK} which states
\begin{lemma2}
Let $A(s) = \sum_{n\leq N} \frac{a(n)}{n^s}$ be a Dirichlet polynomial.  Let $T \geq 1$ and suppose $s_r = \sigma + it_r$ ($r = 1,..., R$) be points with $T<t_1<t_2<...<t_R<2T$ and $t_{r+1} - t_r \geq 1$.  Then
$$\sum_{r=1}^R |A(s_r)|^{2k} \ll \left(T+N^k\right) \log (2N^k)\left(\sum_{n\leq N}d_k(n)|a(n)|^2n^{-2\sigma}\right)^k,
$$where $d_k(m)$ denotes the number of writing $m$ as a product of $k$ numbers, and the implied constant is absolute.
\end{lemma2}
We apply this Lemma to $A(s) = \sum_{n<x^2}\frac{\Lambda_x(n)}{n^s}$, so that $N = x^2 = (\log T)^{4(1+\delta)/\lambda} \ll \log^{8/\lambda} T$.  In our case, $\sigma = 1/2 + \lambda$ so $\sum_{n}d_k(n)|a(n)|^2n^{-2\sigma}$ converges.  Hence
$$\sum_{r=1}^R |A(s_r)|^{2k} \ll_{k, \lambda} T\log T,
$$for any $k$ satisfying $N^k \ll T$.  Since $N$ is a power of $\log T$, we may pick $k$ to be as large a constant as we like.  Now say that $A(s_r) \gg \log ^{1-\delta} T$ for each $1\leq r\leq R$.  Then $R \ll \frac{T}{\log^B T}$ for any fixed $B>0$.  Here, for a given $B$, the choice of $k$ will depend only on $B$ and $\delta$.  Now let
$$E_2  = \{t\in [T, 2T] : |A(\sigma + it)| > \log ^{1-\delta} T\}.
$$From the discussion above, $E_2$ has at most $\ll \frac{T}{\log^BT}$ points that are separated from each other by at least $1$, and thus has $\ll \frac{T}{\log^{B-1}T}$ zeros of $R(t)$ or $I(t)$.  Clearly, the measure of $E_2$ is $\ll \frac{T}{\log^BT}$.

Let $E = E_1 \cup E_2$ be our exceptional set.  Then in intervals outside of $E$, 
\begin{equation}
\frac{d}{dt}\arg \xi(\sigma + it) = \frac{1}{2}\log \frac{t}{2\pi} + O\left(\log ^{1-\delta}T\right),
\end{equation}
where the main term comes from logarithmic differentiation of the gamma factor and the error term from $\zeta '/\zeta(s)$.  Specifically, the main term arises from 
\begin{eqnarray*}
\Im \frac{ds}{dt}\frac{d}{ds}(\log(\pi^{-s/2}) + \log \Gamma(s/2) + s(s-1)/2)
&=& \Im i (-1/2\log \pi +1/2 \log s/2 + O(1/s))\\
&=& 1/2\log(t/2) + O(1).
\end{eqnarray*}
Thus, the argument is a strictly increasing function of $t$ for large $T$ and so the zeros of $R(t)$ and $I(t)$ are simple and strictly interlace.

It remains to examine the normalized spacings.  By integrating (2), we have that for any $\Delta \ll \frac{1}{\log T}$ that
\begin{equation}
\arg \xi(s+i\Delta) - \arg \xi(s) 
= \Delta/2\log(t/2\pi ) + O(\log^{-\delta}T).
\end{equation}
In particular, if $\gamma_1 < \gamma_2$ were consecutive zeros of $R(t)$ in $[T, 2T]$, we claim that $\gamma_2 - \gamma_1 = \frac{2\pi}{\log \frac{\gamma_1}{2\pi}} + O\left(\frac{1}{\log ^{1+\delta}T}\right)$.  Indeed, by the monotonicity of $\arg \xi(\sigma + it)$, 
$$\arg \xi(\sigma+i\gamma_2) - \arg \xi(\sigma + i\gamma_1) = \pi$$.  The claim then follows from (3).  Thus the normalized spacing is $1+O\left(\frac{1}{\log ^{\delta}T}\right)$. Let $\psi$ be the unique zero of $I(t)$ in $(\gamma_1, \gamma_2)$.  Then similarly $\psi - \gamma_1 = \frac{\pi}{\log \frac{\gamma_1}{2\pi}} + O\left(\frac{1}{\log ^{1+\delta}T}\right)$ so that $\psi$ is actually in the middle of $(\gamma_1, \gamma_2)$.  Of course, the same holds with the roles of $R(t)$ and $I(t)$ reversed.  Note that the results above all hold for any $\delta \in (0, 1)$.

\begin{remark}
Our method of proof also provides a direct alternative to the proof of results in \cite{Lag}.  Indeed, if we assume $RH$, then Titchmarch \cite{Titch} in 14.5.1 gives that $\frac{\zeta '}{\zeta}(s) \ll \log ^{1-2\lambda}T$ for $0<\lambda < 1/2$, so if we take $\delta = 2\lambda$, we would have no exceptional set.  For $\lambda > 1/2$, things are even simpler as we have unconditionally that $\frac{\zeta '}{\zeta}(s) \ll 1$.  Finally, for $\lambda = 1/2$, an unconditional result is still possible by using that $\frac{\zeta '}{\zeta}(s) \ll \frac{\log T}{\log \log T}$ for $T\geq 2$, from \cite{Titch} 5.17.4.  Here the error terms with $\frac{1}{\log ^\delta T}$ would be replaced with $\frac{1}{\log \log T}$.  
\end{remark} 

\section{An analogue of Ki's Result}
Here we will provide a sketch of the proof of Part 1 of Theorem 1, leaving out those details which are readily available in \cite{Ki}.  Let 
$$G(z) := \xi(1/2 + iz).
$$Let 
$$H(z) = G(z - i\lambda) + G(z+i\lambda)
$$and
$$J(z) = G(z - i\lambda) - G(z+i\lambda)
$$
By the definition of $G$ and by the functional equation for $\xi$, the two functions above correspond to the real and imaginery parts of $\xi(1/2+\lambda + it)$ when $z = t$ is real.  

The main idea of the proof is to show that there is some small exceptional set, outside which $H(z)$ and $J(z)$ has only real zeros, and so outside of a small exceptional set, $A(s)$ and $B(s)$ only have zeros on the $1/2$ line.  The key to showing this will be to show that $G(z - i\lambda)$ and $G(z+i\lambda)$ have different sizes when $z$ is not real.  The exceptional set will essentially consist of neighbourhoods around the zeros of $\zeta(s)$ for $\Re(s) \geq 1/2+\lambda$.  To be precise, let $\{s_i\}_{i=1}^\infty$ denote those zeros of $\zeta$ with real part $\geq 1/2+\lambda$ and $\Im(s_i) \leq \Im(s_{i+1})$ for all $i$.  Let
$$\R(\kappa) = \{z: |Re(z)| \geq |\Im(s_1)|, \textup{ and } |\Re(z) - \Im(s_k^{*})| > \kappa, s_k^* = s_k \textup{ or } \bar{s_k} \forall k\}.
$$  The exceptional set will be the complement of $\R(\kappa)$.  Our unconditional analogue of Theorem B of \cite{Ki} is below.

\begin{lemma3}
There is some $\kappa>0$ such that all the zeros $z$ with $|z| > \kappa$ of $H(z)$ or $J(z)$ in $\R(\kappa)$ are real.
\end{lemma3}

From now on, we will restrict the discussion to $H(z)$, the situation with $J(z)$ being the same.  The points which are not in $\R(\kappa)$ are all within a constant distance of a zero of $\zeta(s)$ with $\Re(s) \geq 1/2+\lambda$.  The number of such zeros is $\ll T^{1-a\lambda}\log T$ for any $a<1$.  The zeros of $H(z)$ for $T\leq \Re(z) \leq T+C$ is $\ll \log T$ for any constant $C$.  As we mentioned in the Introdution, the proof of this proceeds along standard lines using the argument principle, and we refer the reader to the proof of Proposition 2.2 of \cite{Ki} for more details.  Thus, the number of zeros of $H(z)$ outside of $R(\kappa)$ is $\ll T^{1-a\lambda}\log^2T$, uniformly in $\lambda$.  

The proof of Lemma 3 is the same as the proof of Theorem B in \cite{Ki}.  We provide a sketch below, referring the reader to details which are available in \cite{Ki} whenever possible.

\begin{proof}
Set $z = x-iy \in R(\kappa)$ with $x>1$ and $y>0$.
We have by Equation 3.6 and 3.7 in \cite{Ki} that for some constants $H_0, H_1, \sigma_1>0$, explicitly described in \cite{Ki}, 
\begin{eqnarray*}
\log \left| \frac{G(\bar{z} - i\lambda)}{G(z - i\lambda)}\right|
&\leq& 4y\left( H_0 \frac{\log(\kappa + x)}{\kappa^2} + H_0 \sum_{n\geq 1} \frac{\log x}{(\kappa + n)^2} - \frac{\lambda H_1 \log x}{1+\sigma_1^2}\right)\\
&=& 4y\log x \left( H_0 \left(\frac{1}{\kappa} + \frac{1}{\kappa^2}\right) - \frac{\lambda H_1}{1+\sigma_1^2}\right) (1+o(1))\\
&=:& 4y\log x \delta  (1+o(1))
\end{eqnarray*}
where we have set $\delta = H_0 \left(\frac{1}{\kappa} + \frac{1}{\kappa^2}\right) - \frac{\lambda H_1}{1+\sigma_1^2}$.  We then pick $\kappa$ such that $\delta> 0$ and so
$$\left|\frac{G(\bar{z} - i\lambda)}{G(z - i\lambda)}\right|
\leq e^{-4y\log x\delta}.
$$
Now we have that
\begin{eqnarray*}
|H(z)| &=& G(z - i\lambda) + G(z+i\lambda)\\
&\geq& |G(z - i\lambda)|\left(1 -  \frac{G(z+i\lambda)}{G(z - i\lambda)}\right)\\
&=& |G(z - i\lambda)|\left(1 -  \frac{G(\bar{z}-i\lambda)}{G(z - i\lambda)}\right)\\
&>& 0,
\end{eqnarray*}for $z \in \R(\kappa)$, and for $|z| \geq \kappa$.
Thus $H(z) \neq 0$ for $z \in \R(\kappa)$, $|z| \geq \kappa$ and $y>0$.  Since $H(z) = \overline{H(\bar{z})}$ and $H(-z) = H(z)$, we have that the same result holds true for $y<0$, from which we deduce that all zeros of $H(z)$ with for $z \in \R(\kappa)$ and $|z| \geq \kappa$ must be real
\end{proof}

\end{document}